\documentclass[11pt]{amsart}
\usepackage{graphicx}

\usepackage{amsmath}
\usepackage{amsthm}
\usepackage{amssymb}
\usepackage{enumerate}
\usepackage{epic,eepic}

\newtheorem*{thma}{Theorem A}
\newtheorem*{thmb}{Theorem B}

\def\Z{\mathbb Z}
\def\R{\mathbb R}
\def\Q{\mathbb Q}
\def\N{\mathbb N}
\def\A{\mathcal A}

\def\C{\mathcal C}
\def\L{\mathcal L}
\newcommand{\mat}[1]{\boldsymbol{#1}}
\newtheorem{lem}{Lemma}[section]
\newtheorem{notat}{Notation}[section]
\newtheorem{prop}[lem]{Proposition}
\newtheorem{coro}[lem]{Corollary}

\newtheorem{de}[lem]{Definition}

\newtheorem{pozn}[lem]{Remark}

\def\pf{\begin{proof}}
\def\pfk{\end{proof}}

\begin{document}
\title{Sturm numbers and substitution\\ invariance of 3iet words}

\author[Arnoux]{ Pierre Arnoux}

\address[Arnoux]{ Institut de Math\'ematiques de Luminy, CNRS UPR 9016\\
163, avenue de Luminy, Case 907,
13288 Marseille cedex 09
France}
\email[Arnoux]{ arnoux@iml.univ-mrs.fr}

\author[Berth\'e]{Val\'erie~Berth\'e}

\address[Berth\'e]{Laboratoire d'Informatique, de Robotique et de MicroŽlectronique de Montpellier
161, rue Ada, 34 392 Montpellier Cedex 5, France}
\email[Berth\'e]{berthe@lirmm.fr}

\author[Mas\'akov\'a]{ Zuzana~Mas\'akov\'a}
\address[Mas\'akov\'a]{ Doppler Institute %%for Mathematical Physics and Applied Mathematics\\
\& Department of Mathematics, FNSPE Czech Technical University, Trojanova 13, 120 00 Praha 2, Czech Republic}
\email[Mas\'akov\'a] {masakova@km1.fjfi.cvut.cz}

\author[Pelantov\'a]{Edita~Pelantov\'a}
\address[Pelantov\'a]{ Doppler Institute %%for Mathematical Physics and Applied Mathematics\\
\& Department of Mathematics, FNSPE Czech Technical University, Trojanova 13, 120 00 Praha 2, Czech Republic}
\email[Pelantov\'a] {pelantova@km1.fjfi.cvut.cz}

\date{\today}

\begin{abstract} In this paper, we give a necessary condition for an infinite word defined by a
non-degenerate  interval exchange on three intervals  (3iet word) to be invariant by a substitution: a natural  parameter associated to this word must be a Sturm number. We deduce some
algebraic  consequences from this condition concerning the incidence matrix
of the associated  substitution. As a by-product of our proof, we give a combinatorial characterization of 3iet words.
\end{abstract}
\maketitle

\section{Introduction}

The original definition of a Sturm number using continued fractions was introduced in 1993
when Crisp et al.~\cite{crisp} showed that a homogeneous sturmian word (i.e.,  a sturmian
word with slope $\varepsilon$ and intercept $x_0=0$) is invariant under a non-trivial
substitution if and only if $\varepsilon$ is a Sturm number. In 1998, Allauzen~\cite{allauzen}
provided a simple characterization of Sturm numbers:

{\em A quadratic irrational number $\varepsilon$ with conjugate $\varepsilon'$ is called
a Sturm number if
\begin{equation}\label{e:000}
\varepsilon\in(0,1) \quad\hbox{ and }\quad\varepsilon'\notin(0,1)\,.
\end{equation}
}

For general sturmian words (with arbitrary intercept $x_0$), the fact that $\varepsilon$ is
a Sturm number is only a necessary but not a sufficient condition for invariance under a substitution; this is clear since there can be only a countable number of such invariant words, while the sturmian words with a given slope are determined by their intercept, hence they are uncountable in number. For a complete characterization, see \cite{Yasutomi,BaMaPe,berthe}.

In this paper we study invariance under substitution of infinite words coding non-degenerate
exchange of three intervals with permutation (321).\footnote{An exchange of intervals
is non-degenerate if it satisfies i.d.o.c. \cite{kean}. For more details, see
Section \ref{sec:subsec4.3}} These words, which are here called
non-degenerate 3iet words, are one of the possible generalizations of sturmian words to a
three-letter alphabet.
Some  combinatorial properties characterizing the language of 3iet words are described in~\cite{FeHoZa}.  It is  well  known  that   substitutive   3iet words,
 that is,  3iet  words that are image by a morphism of a  fixed point of a substitution
 correspond to  quadratic parameters, see e.g.  \cite{adamczewski,AS,BC,FeHoZa,PLV}.
Let us stress the fact  that we  consider  in the present paper  fixed points  of substitutions
and not  substitutive words.

Sturmian words can be equivalently defined as aperiodic words coding a rotation, that is, an exchange of two intervals with lengths say $\alpha$, $\beta$. The slope of the sturmian word, which we have denoted by
$\varepsilon$, is then equal to $\varepsilon=\frac{\alpha}{\alpha+\beta}$. The term `slope' for
the parameter $\varepsilon$ comes from the fact that the sturmian word with slope $\varepsilon$
can be constructed by projection of points of the lattice $\Z^2$ to the straight line
$y=\varepsilon x$; it will prove convenient to abuse the language by speaking of the {\em slope of the rotation}, this slope is the complement to 1 of the more usual {\em angle} of the rotation. Since the sturmian word does not depend on the absolute lengths of the
two intervals being exchanged but one their ratio,  then the lengths are often normalized to satisfy $\alpha+\beta=1$.
In this case $\alpha$ and $\varepsilon$ coincide.

The same situation appears for 3iet words which code exchange of three intervals with lengths, say,
$\alpha$, $\beta$, $\gamma$. The commonly used normalization of parameters is
$\alpha+\beta+\gamma=1$. However, much more suitable appears to be the normalization
$\alpha+2\beta+\gamma=1$. Let us mention three arguments in favor as results of
papers~\cite{adamczewski,3iet,FeHoZa,GuMaPe}. If $u$ is an infinite word coding
exchange of three intervals of lengths $\alpha$, $\beta$, $\gamma$, then:

\begin{itemize}
\item the infinite word $u$ is aperiodic
if and only if $\displaystyle{\frac{\alpha+\beta}{\alpha+2\beta+\gamma} \ \notin \ \Q}$;

\item  if $u$  is assumed to be aperiodic, then $u$
codes a non-degenerate exchange of three intervals
if and only if $\displaystyle{\frac{\alpha+\beta+\gamma}{\alpha+2\beta+\gamma}} \ \notin \
\Z+\Z\displaystyle{\frac{\alpha+\beta}{\alpha+2\beta+\gamma}}$;

\item
the infinite word $u$ can be constructed by projection of points of the lattice $\Z^2$
on the straight line $y=\displaystyle{\frac{\alpha+\beta}{\alpha+2\beta+\gamma}\ x}$.

\end{itemize}
 We will give  in Section \ref{sec:sec3} and \ref{sec:sec4}
 a short proof of these facts by recalling that such an exchange of three intervals can always be obtained as an induced map of a rotation (exchange of two intervals) on an interval of length $\alpha+2\beta+\gamma$; in the process, we will give a complete combinatorial characterization of 3iet words, as follows:

 \begin{thma} Let $u$ be a sequence on the alphabet $\{A,B,C\}$ whose letters  have positive densities. Let
 $\sigma: \{A,B,C\}^*\to \{0,1\}^*$$\sigma$ and $\sigma: \{A,B,C\}^*\to \{0,1\}^*$
  be the morphisms  defined by $$\sigma(A)=0,\quad \sigma(B)= 01,\quad \sigma(C)=1$$  $$\sigma'(A)=0,\quad\sigma'(B)= 10,\quad \sigma'(C)=1.$$

 The sequence $u$ is an aperiodic 3iet word if and only if $\sigma(u)$ and $\sigma'(u)$ are sturmian words.
 \end{thma}

This paper adds yet another argument supporting the normalization $\alpha+2\beta+\gamma=1$, by the following necessary condition:

\begin{thmb}
If a non-degenerate 3iet word is invariant under a primitive substitution, then
$$
\varepsilon:=\displaystyle{\frac{\alpha+\beta}{\alpha+2\beta+\gamma}} \quad
\hbox{is a Sturm number.}
$$
\end{thmb}

Remark that, in that case, the corresponding homogeneous sturmian
word is also substitution invariant. A forthcoming paper  \cite{BaMaPe2}
 will give a complete characterization of
substitution invariant 3iet words. Note that it is a natural question to ask wether, when $u$ is a substitution
invariant 3iet word, one or both of the sturmian words
$\sigma(u)$, $\sigma'(u)$ are also substitution invariant.

This paper is organized as follows. The introductory notation and
definitions  are  given in Section \ref{sec:sec2}.  Section
\ref{sec:sec3} and \ref{sec:sec4} are devoted to the  description
of  a classical  exduction process  in terms of  substitutions and
to the proof of Theorem A. Section \ref{sec:sec5} and Section
\ref{sec:sec6} gather the  required material  for the proof of
Theorem  B, namely, properties  of translation vectors and balance
properties. Theorem B is proved in Section \ref{sec:sec7}.

%%%%%%%%%%%%%%%%%%%%%%%%%%%%%%%%%%%%%%%%%%%%%%%%%%%%%%%%%%%%%%%%%%%%%%%%%%%
\section{Preliminaries} \label{sec:sec2}

We work with finite and infinite words over a finite
alphabet $\A=\{a_1,\dots,a_k\}$. The set of
all finite words over $\A$ is denoted by $\A^*$. Equipped with the binary operation
of concatenation and the empty word, it is a free
monoid. The length of a word
$w=w_1w_2\cdots w_n$ is denoted by $|w|=n$, the number of letters
$a_i$ in the word $w$ is denoted by $|w|_{a_i}$.

An infinite concatenation of letters of $\A$ forms the infinite
word $u=(u_n)_{n\in\N}$,
$$
u = u_0u_1u_2\cdots\,.
$$

A word $w$  is said to be a \emph{factor} of a word
$u=(u_n)_{n\in\N}$ if there is an index $i\in\N$ such that $w=u_i
u_{i+1}\cdots u_{i+n-1}$. The set of all factors of $u$ of length
$n$ is denoted by $\mathcal{L}_n(u)$. The \emph{language}
$\mathcal{L}(u)$ of an infinite word $u$ is the set of all its
factors, that is,
$$
\mathcal{L}(u) = \bigcup_{n\in\N} \mathcal{L}_n(u)\,.
$$

The \emph{(factor) complexity} $\C_u$ of an infinite word $u$ is the function
$\C_u:\N\rightarrow\N$ defined as
$$
\C_u(n) := \#\mathcal{L}_n(u)\,.
$$
The density of a letter $a\in\A$,
representing the frequency of occurrence of the letter $a$ in an
infinite word $u$, is defined by
$$
\rho(a) := \lim_{n\rightarrow\infty} \frac{\#\{i\ |\ 0\leq i< n,\
u_i=a\}}{n}\,,
$$
if the limit exists (this is always the case for 3iet words, as is easy to prove).

 Let ${\mathcal A}$ and ${\mathcal B}$ be two alphabets.
 A mapping $\varphi:\A^*\rightarrow  {\mathcal B}^*$ is said to be a
\emph{morphism}  if
$\varphi(w\widehat{w})=\varphi(w)\varphi(\widehat{w})$ holds for
any pair of finite words $w,\widehat{w}\in\A^*$. Obviously, a
morphism is uniquely determined by the images $\varphi(a)$ for all
letters $a\in\A$. If ${\mathcal A}$ and ${\mathcal B}$ coincide and if the images
of    the letters   are  never  equal to the empty word,
then $\varphi$ is called a \emph{substitution}.

The action of a morphism $\varphi$ can be naturally extended to
infinite words by the prescription
$$
\varphi(u) = \varphi(u_0u_1u_2\cdots) :=
\varphi(u_{0})\varphi(u_{1})\varphi(u_2)\cdots\,.
$$
A infinite word $u\in\A^\N$ is said to be a \emph{fixed point} of
the morphism $\varphi$ if $\varphi(u)=u$.

The incidence matrix of a morphism $\varphi$ over the
alphabet $\A$ is an important tool which brings a lot of information about the
combinatorial properties of the fixed points of the morphism. It is defined by
$$
(\mat{M}_\varphi)_{ij} = |\varphi(a_i)|_{a_j} =
\text{number of letters $a_j$ in the word $\varphi(a_i)$}\,.
$$
A morphism $\varphi$ is called primitive if there exists an integer $k$ such that
the matrix $\mat{M}_{\varphi}^k$ is positive.

Assume that an
infinite word $u$ over the alphabet $\A=\{a_1,\dots, a_k\}$
is a fixed point of a primitive substitution $\varphi$. It is known~\cite{Queffelec} that
in such a case the densities of letters in $u$  are well defined.
The vector
$$
\vec{\rho}_u = \bigl(\rho(a_1),\dots,\rho(a_k)\bigr)\,.
$$
is a left eigenvector of the incidence matrix $\mat{M}_\varphi$,
i.e.,  $\vec{\rho}_u\mat{M}_\varphi=\Lambda\vec{\rho}_u$. Since the
incidence matrix $\mat{M}_\varphi$ is a non-negative integral
matrix, we can use the Perron-Frobenius Theorem stating that
$\Lambda$ is the dominant eigenvalue of $\mat{M}_\varphi$.
Moreover, all eigenvalues of $\mat{M}_\varphi$ are algebraic
integers.

%The right eigenvector of the incidence matrix corresponding to the
%dominant eigenvalue plays an
%important role for the \emph{geometric representation} of a fixed
%point of a morphism. Let $\vec{x}=(x_1,x_2,\ldots,x_k)^T$ be a positive
%right eigenvector of $\mat{M}_\varphi$
%corresponding to $\Lambda$.
%
%If the biinfinite word $u = \cdots
%u_{-3}u_{-2}u_{-1}|u_0u_1u_2\cdots$ is the fixed point of $\varphi$, we can define
%a discrete set $\Sigma$ by
%$$
%\Sigma = \{t_n\ |\ n\in\Z\}\,,\quad \text{where} \quad t_0=0 \
%\hbox{ and } \ t_{n+1} - t_n = x_i\ \Leftrightarrow\ u_n=a_i\,.
%$$
%Since $u$ is a fixed point of a morphism, the construction of
%$\Sigma$ implies that $\Lambda\Sigma \subset \Sigma$, i.e.\ $\Sigma$ is a self-similar set.

%%%%%%%%%%%%%%%%%%%%%%%%%%%%%%%%%%%%%%%%%%%%%%%%%%%%%%%%%%%%%%%%%%%%%%%%%%%%%%%%%%%%%%%%
\section{Exchanges of three intervals as induction of rotations} \label{sec:sec3}
Let $\alpha,\beta,\gamma>0$ and denote by
$$
I_A:=[0,\alpha), \quad I_B:=[\alpha,\alpha+\beta), \quad
I_C:=[\alpha+\beta,\alpha+\beta+\gamma),\ \hbox{ and }\ I:=I_A\cup
I_B\cup I_C\,,
$$
and let $t_A=\beta+\gamma, t_B=\gamma-\alpha , t_C=-\alpha-\beta \in\R$ be translations vectors; we have:
$$
I_A\cup I_B\cup I_C = (I_A+t_A)\cup (I_B+t_B)\cup (I_C+t_C)\,.
$$

The map $T$ defined on $I$ by $T(x)=x+t_X$ if $x\in I_X$, $X=A,B,C$ is the exchange of
three intervals $I_A$, $I_B$, $I_C$ with the permutation (321).

As was already known long time ago  (see \cite{KASTE}), this map can be obtained as the induction  of a rotation on a suitable interval. We recall the construction; let $I_D=[\alpha+\beta+\gamma,\alpha+2\beta+\gamma)$, and define $J=I\cup I_D$. Let $R$ be the rotation of angle $\frac {\beta+\gamma}{\alpha+2\beta+\gamma}$ on $J$, defined by $R(x)=x+\beta+\gamma$ if $x\in I_A\cup I_B$, and $R(x)=x-\alpha-\beta$ if $x\in I_C\cup I_D$; it exchanges the two intervals  $I_A\cup I_B$ and $I_C\cup I_D$.

For a subset $E$ of $X$, the {\em first return time} $r_{E}(x)$ of
a point $x\in E$ is defined as  $\min\{n>0\mid R^nx\in E\}$. If
the return time is always finite, we define the {\em induced map}
or {\em first return} map  of $R$ on $E$ by
$R_{E}(x)=R^{r_{E}(x)}(x)$.
\begin{lem} \label{lem:lem3.1}
 The map $T$ is the first return map of $R$ on $I$.
\end{lem}
\pf
 One checks that $R(I_A)=[\gamma+\beta, \alpha+\gamma+\beta)=T(I_A)$,  $R(I_C)=[0, \alpha)=T(I_C)$, $R(I_B)=I_D$, and $R^2(I_B)=R(I_D)=T(I_B)$.
\pfk

Hence $T$ can be obtained as induction of a rotation on a left
interval (for more details, see e.g. \cite{adamczewski,FeHoZa} or the survey \cite{BFZ}).
It  can also be obtained as induction on a right
interval, and this remark will prove important below: define
$I_E=[-\beta, 0)$, and $J'=I_E\cup I$; consider the rotation $R'$
on $J'$ by the same angle
$\frac{\beta+\gamma}{\alpha+2\beta+\gamma}$; in the same way, one
proves that $T$ is obtained as the first return map of $R'$ on
$I$.

The underlying rotation $R$ turns out to play an important role in
the study of $T$; this explains the appearance of the number
$\varepsilon=\frac{\alpha+\beta}{\alpha+2\beta+\gamma}$ in the
introduction: it is the slope of the rotation  $R$.

\begin{notat} From now on, we will take the normalization
$\alpha+2\beta+\gamma=1$.
\end{notat}

This amounts to normalize the interval of definition of $R$ to 1, and will greatly simplify the notation below.

%%%%%%%%%%%%%%%%%%%%%%%%%%%%%%%%%%%%%
\section{Characterization of non-degenerate 3iet words}\label{sec:sec4}

\subsection{From 3iet words to sturmian words}
With an initial point $x_0\in I$ we associate an infinite word
which codes the orbit of $x_0$ under $T$ with respect to the
natural partition in three intervals (see Definition \ref{def:def4.1}  below). It turns out to be useful to
shift the interval of definition, so that the free choice of the
initial point $x_0$ for the orbit is replaced by the choice of  a parameter that we call $c$
as the position of the interval. The initial point for the orbit
thus can always be chosen as the origin. For this we introduce the
new parameters
$$
\varepsilon := \alpha+\beta\,,\qquad
l:= \alpha+\beta+\gamma\,,\qquad
c:= -x_0\,,
$$
The number $\varepsilon$ is the slope of the underlying rotation
$R$, and $l$ determines the length of the induction interval $J$.
It is obvious that the above parameters satisfy
\begin{equation}\label{e:100}
\varepsilon\in(0,1)\,,\quad \max(\varepsilon,1-\varepsilon) < l<1\,,\quad
-l<c\leq 0 .
\end{equation}

We redefine in this setting five intervals $$I_A=[c, c+\alpha), \
I_B=[c+\alpha, c+\varepsilon), \ I_C=[c+ \varepsilon,c+l), \
I_D=[c+l, c+1), \ I_E=[c-\beta,c).$$ We define $I=I_A\cup I_B\cup
I_C$; the map $T$ (introduced   above  in Section \ref{sec:sec3}) is defined on $I$ as the exchange of three
intervals $I_A, I_B, I_C$ according to the permutation (321).

We also define $$J_0=I_A\cup I_B, \ J_1=I_C\cup I_D \mbox{   and } J=J_0\cup J_1$$
$$J'_0=I_E\cup I_A, \ J'_1=I_B\cup I_C \mbox{ and }
J'=J'_0\cup J'_1.$$  The rotation $R$ (resp. $R'$)  is then defined on $J$  (resp. $J'$) by
the exchange of $J_0$ and $J_1$ (resp.  $J'_0$ and $J'_1$); it has  angle
$1-\varepsilon$,   $J_0$ and $J'_0$ have length  $ \varepsilon$,
whereas  $J_1$ and $J'_1$ have length $1-\varepsilon$.

Let us formulate the definition of 3iet words with the use of  these new parameters.

\begin{de}\label{def:def4.1}
Let $\varepsilon,l,c\in\R$ satisfy~\eqref{e:100}.
The infinite word $(u_n)_{n\in\N}$ defined by
\begin{equation}\label{e:103}
u_n=\left\{\begin{array}{cl}
A& \hbox{if } \ T^n(0)\in I_A\,,\\
B& \hbox{if } \ T^n(0)\in I_B\,,\\
C& \hbox{if } \ T^n(0)\in I_C
\end{array}\right.
\end{equation}
is called the 3iet word with parameters $\varepsilon,l,c$.
\end{de}

There is  classical  and simple way to give a
combinatorial interpretation of   the induction process of
Lemma \ref{lem:lem3.1}  in terms of substitutions.
Consider indeed the orbit $(T^n(0))_{n\in\N}$ of 0 under $T$; it is
clear by Lemma \ref{lem:lem3.1}  that it is a subset of the orbit $(R^n(0))_{n\in\N}$ under
$R$; the points of the second orbit which are not in the first are
exactly the points in $I_D$, and their preimages are exactly the
points in $I_B$; the return time of these points to $I$ is 2. Let
$u$ be the coding of the orbit of 0 under $T$ with respect to the
partition in three intervals $I_A, I_B,I_C$; to obtain the coding
of the orbit of the same point under $R$, with respect to the
partition $I_A, I_B,I_C,I_D$, this argument shows that it is
enough to introduce a letter $D$ after each $B$, that is to
replace $B$ by $BD$; to obtain the natural sturmian coding with
respect to the partition $J_0, J_1$, we then project letters $A,B$
to 0 and $C,D$ to 1.

\begin{de}
We denote by $\sigma$ (resp. $\sigma'$) the morphism from
$\{A,B,C\}^*$ to $\{0,1\}^*$ defined by $\sigma(A)=0$,
$\sigma(B)=01$, $\sigma(C)=1$ (resp. $\sigma'(A)=0$,
$\sigma'(B)=10$, $\sigma'(C)=1$).
\end{de}

We  thus have proved the following:

\begin{lem}\label{l:nove}
Let $u$ be the coding of the orbit of 0 under $T$, with respect to
the partition $I_A, I_B,I_C$, and let $v$ (resp. $v'$) be the
coding of the orbit of 0 under $R$ (resp. $R'$) with respect to
the partition $J_0,J_1$ (resp. $J'_0,J'_1$). Then $v=\sigma(u)$,
$v'=\sigma(u')$.

This implies that, if $\varepsilon$ is irrational, then   $\sigma(u)$
and $\sigma(u')$ are sturmian sequences whose density of $0$ equals
$\varepsilon$.
\end{lem}

\subsection{Characterization theorem} \label{sec:sec4.2}
We will now prove the reciprocal (Theorem A below);  we need some properties of
sturmian sequences.

Let $v$ be a sturmian sequence  that codes the orbit of
  a rotation of angle $1-\varepsilon$  modulo $1$
  with   density of $0$  equal to
$\varepsilon$, and let $V_n$ be the prefix of $v$ of length $n$, i.e.,
$V_n =v_0 \cdots v_{n-1}$.
Define a map $$f:\{0,1\}^*\to \R \mbox{ by }
f(V)=|V|_0(1-\varepsilon)-|V|_1\varepsilon .$$ From the definition,
we see that
 $$\forall n, \  f(V_n)=|V_n|_0 - n \varepsilon,$$ hence   the sequence $(f(V_n))_{n\in\N}$ is the orbit of 0
under a rotation defined on an interval $[c,c+1)$, with $c=\inf\{ f(V_n)| n\in\N\}$; in particular,
we have, for all integers $i,j$,  $|f(V_i)-f(V_j)|<1$.

We have the following lemma:

\begin{lem}\label{lem:lem4.4}
Let $v$ be a sturmian sequence, and let $(n_k)_{k\in\N} $ be a
strictly increasing sequence of integers  that satisfies $v_{n_k-1}=0$,
$v_{n_k}=1$. Define a new sequence $v'$ by: $v'_{n_k-1}=1$,
$v'_{n_k}=0$, $v'_i=v_i$ otherwise. The sequence $v'$ is sturmian
if and only if for every $i$ which is not in the sequence
$(n_k)_{k\in\N}$, and for all $j$, we have  $f(V_{n_j})>f(V_i)$.
\end{lem}

\pf   For all $n$, let $V'_n$  stand for  the prefix of $v'$  of length $n$.
We have
$f(V_i)=f(V'_i)$, except if $i=n_k$, in which case one checks that
$f(V'_{n_k})=f(V_{n_k})-1$.

 We  first assume that $v'$ is sturmian.
 Suppose that $f(V_i) \geq  f(V_{n_j})$, for some $i,j$,  with $i$ not in the sequence
$(n_k)_{k\in\N}$; then we must
have $f(V'_i) \geq f(V'_{n_j})+1$; but this is impossible since  $v'$ is a
sturmian sequence with  same density of $0$'s as $v$. Hence  for every $i$ which is not in the sequence
$(n_k)_{k\in\N}$, and for all $j$, we have $f(V_i)  <  f(V_{n_j})$.

Conversely, we assume that  for every $i$ which is not in the sequence
$(n_k)_{k\in\N}$, and for all $j$, we have $f(V_{n_j})>f(V_i)$. One checks that  for all integers $i,j$,  $|f(V'_i)-f(V'_j)|<1$. Indeed, this is immediate
if  $i$ and $j$  belong simultaneously  to  $(n_k)_{k\in\N}$, or else if
none of them belongs to  this sequence.   If $i$  is not in the sequence
$(n_k)_{k\in\N}$, then $  |f(V'_i)-f(V'_{n_j})|= | f(V_i) -f(V_{n_j}) +1 |   < 1$.
We deduce  that  the sequence $v'$ is a  balanced sequence.
Indeed,   take two factors $W$ and $W'$  of the same length   $n$ of  the sequence
$v'$  that occur  respectively at index
$i$ and $j$.
One has  $$ ||W|_0 -|W'|_0 |
=|(f(V'_{i+n}) - f( V'_{j+n} ))- (f(V'_{i})-f(V'_j))|  < 2.$$
We deduce that  the densities of letters  are well-defined in $v'$. By construction, they
coincide with the   densities of letters  for the sequence $v$, hence $v'$ is  an  aperiodic
balanced sequence,
 it is thus  a sturmian sequence, according to~\cite{Hedlund}. \pfk

We are now in position the prove the first theorem:

 \begin{thma}
 Let $u$ be a sequence on the alphabet $\{A,B,C\}$ whose letters  have positive  densities.
 This sequence is an aperiodic 3iet word if and only if $\sigma(u)$ and $\sigma'(u)$ are sturmian words.
 \end{thma}
\pf We have proved above (Lemma~\ref{l:nove}) that the condition
is necessary. Let us prove it is sufficient. Let $v=\sigma(u),
v'=\sigma'(u)$ be the two sturmian words; by construction, they
have the same slope $\varepsilon$, and they coincide except on a
sequence of pairs of indices $(n_k-1, n_k)$, corresponding to the
images of $B$, where 0 is replaced by 1 and vice versa.

Define the function $f$ as above, and define $c=\inf \{f(V_k)\mid
k\in\N\}$, $l=\inf \{f(V_{n_k})\mid k\in\N\}-c$. From Lemma
\ref{lem:lem4.4}, we deduce that an index $i$ is of the form $n_k$
if and only if $f(V_i)   \geq  l+c$ if $ \inf \{ f(V_{n_k}) \} =
\min \{ f(V_{n_k}) \} $ (resp.  $f(V_i)   >   l+c$ otherwise).
Then,
 one  checks that the sequence $u$ is generated by the exchange $T$  of the
three intervals $I_A, I_B, I_C$
with  either
$$I_A=[c,c+\varepsilon+l-1), \ I_B=[c+\varepsilon+l-1,
c+\varepsilon), \ I_C=[c+\varepsilon, c+l)$$
or
$$I_A=(c,c+\varepsilon+l-1], \ I_B=(c+\varepsilon+l-1,
c+\varepsilon], \ I_C=(c+\varepsilon, c+l],$$
the   choice of the intervals  being determined   by the values
of  $u$,  and thus of $v$ and $v'$,  at    the indices   (if any)   where
the orbit of $0$ under  $T$  meets   discontinuity points.
 The interval $I_B$ corresponds to
the times $n_k-1$ for the sequence $v$, and $I_A$ and $I_C$ resp.
to value 0 and 1 for the other times.  We deduce that $u$  is  aperiodic  from
the  irrationality  of $\varepsilon$. \pfk

\begin{figure}[ht]
\begin{center}
\includegraphics[height=8cm]{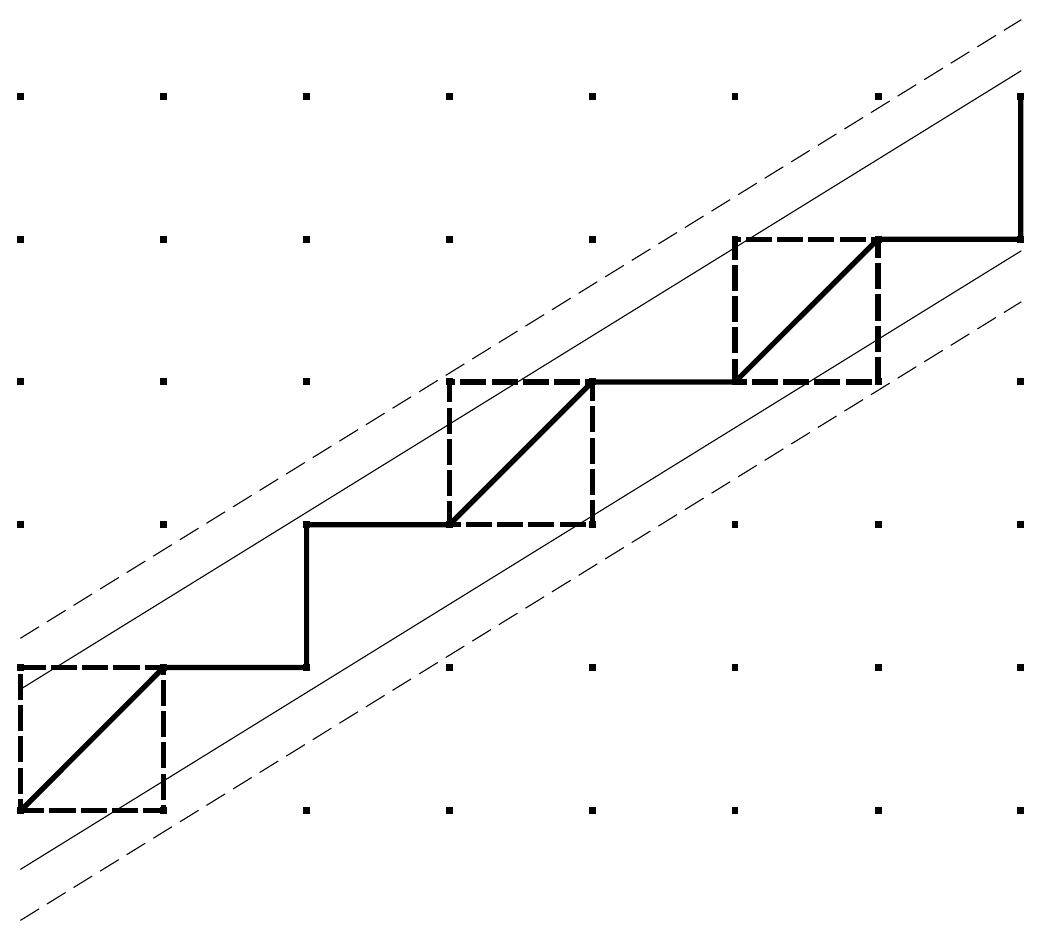}
\caption{The stepped line associated to a 3-iet word and its two sturmian extensions}
\label{fig:line}
\end{center}
\end{figure}

Figure \ref{fig:line} gives a geometric interpretation of the
proof; to the 3iet word $u$, we have associated a stepped line
(bold line), by associating letter $A$ to   vector $(1,0)$, $B$ to $(1,1)$,
and $C$ to $(0,1)$. Remark that this stepped line is contained in
a ``corridor" of width less than 1; in dashed lines are shown the
two sturmian lines associated to $v$ and $v'$, obtained by
enlarging the corridor on the right or the left to the width of
the unit square.

\subsection{Complexity}\label{sec:subsec4.3}

It is known that the factor complexity of the infinite
words~\eqref{e:103} satisfies $\C(n)\leq 2n+1$ for all $n\in\N$. A
short proof can be given by considering the partition $\mathcal P$
in three intervals; to count the number of factors of length $n$,
it is enough to count the number of atoms of the partition
$\bigvee_{k=0}^{n-1}T^{-k}\mathcal P$. But it is easy to prove
that these atoms are intervals, bounded by reciprocal images of
the two discontinuity points. As there can be at most $2n$ such
points between time 0 and $n-1$, there are at most 2n+1 intervals.

The infinite words $(u_n)_{n\in\N}$ which have full complexity are
called {\em non-degenerate} (or regular) 3iet words; 3iet words
for which there exists $n$ such that $\C(n)<2n+1$ are called
degenerate.

The necessary and sufficient condition for a word $(u_n)_{n\in\N}$
coding 3iet to be non-degenerate is the so-called i.d.o.c. (infinite distinct orbit  condition). This
notion has been introduced by Keane~\cite{kean} and requires, in this case,  that
the orbits of the two points of discontinuity of the
transformation $T$ are disjoint, formally
$\{T^{n}(c+l-1+\varepsilon)\}_{n\in\N}\cap
\{T^{n}(c+\varepsilon)\}_{n\in\N}=\emptyset$. If this condition holds true, then the partition above is limited by exactly $2n$ points on the interval, hence has $2n+1$ atoms. The condition
i.d.o.c is equivalent to
\begin{equation}\label{e:101}
\varepsilon\notin\Q\quad\hbox{and}\quad l\notin\Z+\Z\varepsilon=:\Z[\varepsilon]\,,
\end{equation}
see~\cite{adamczewski,GuMaPe}.

\begin{pozn}
If $\varepsilon$ is irrational, it is classical that the rotation
$R$ is uniquely ergodic, which implies that $T$ is also uniquely
ergodic. In that case, the densities of letters in the 3iet
aperiodic word are well defined and $\vec{\varrho}_u$ is
proportional to the vector of lengths of intervals $I_A$, $I_B$,
$I_C$.

If $\varepsilon$ is rational, the sequence $u$ is periodic, hence the densities  exist in a trivial way.
\end{pozn}

\section{Translation vectors}\label{sec:sec5}
Let $u$ be a 3iet word  such as defined in Definition \ref{def:def4.1}. In our considerations, the column vector of translations will play a crucial role. We denote it by
$$
\vec{t}=\left(\!\!\!\begin{array}{c}
t_A\\t_B\\t_C
\end{array}\!\!\!\right)=
 \left(\!\!\!\begin{array}{c}
1-\varepsilon\\1-2\varepsilon\\-\varepsilon
\end{array}\!\!\!\right)\,.
$$

A first remark is that the vector of translations is orthogonal to
the vector of densities; this can be checked directly, and
interpreted as the fact that the mean translation is 0, because
the orbit under the action of the map $T$ is bounded.

We assume furthermore  that
$u$ is fixed  by some substitution $\varphi$.
We will now obtain a more subtle equation, using the substitution $\varphi$.
Let us define a function $g$  (in the  flavour of the map $f$ defined in Section
\ref{sec:sec4.2}) on the prefixes of the infinite word
$u$, the fixed point of $\varphi$. For the prefix $w=u_0u_1\cdots
u_{n-1}$, $n\geq 0$, we put
$$
g(u_0u_1\cdots u_{n-1}):=T^n(0)=|w|_At_A+|w_B|t_B+|w_C|t_C\,.
$$
In particular,  the image of the empty word equals $0$\,. For $X\in\{A,B,C\}$, put
$$
E_{X}:=\Bigl\{ g\bigl(u_0u_1\cdots u_{n-1}\bigr) \Bigm|
u_n=X\Bigr\}=\left\{(|w|_A,|w|_B,|w|_C)\, \vec{t}\,\Bigm|\,
wX\hbox{ is a prefix of } u\right\}\,.
$$
Clearly, the closure of the set $E_X$ satisfies
$\overline{E}_X=I_X$.

The infinite word $u_0u_1u_2\cdots
=\varphi(u_0)\varphi(u_1)\varphi(u_2)\cdots$ can be imagined as a
concatenation of blocks $\varphi(A)$, $\varphi(B)$, $\varphi(C)$.
Positions, where these blocks start, and the corresponding
iterations of $T$,  are given by the following sets. For
$X\in\{A,B,C\}$, put
$$
E_{\varphi(X)}:=\Bigl\{ g\bigl(\varphi(u_0u_1\cdots u_{n-1})\bigr)
\Bigm| u_n=X\Bigr\} = \Bigl\{ g\bigl(\varphi(w)\bigr) \Bigm|
wX\hbox{ is a prefix of } u\Bigr\}\,.
$$
From the definition of the matrix $\mat{M}_\varphi$ it follows that
\begin{equation}\label{eq:nojo}
E_{\varphi(X)}=\left\{(|w|_A,|w|_B,|w|_C) \mat{M}_\varphi
\vec{t}\,\Bigm|\, wX\hbox{ is a prefix of } u\right\}\,.
\end{equation}
Obviously,
$$
E_{\varphi(A)}\cup E_{\varphi(B)} \cup E_{\varphi(C)} \subset \{T^n(0) \mid n\in\N\}\subset I\,,
$$
and the union is disjoint. The fact that $T^{k}(0)$ belongs to $E_{\varphi(A)}$ is equivalent to
\begin{itemize}
\item
$u_ku_{k+1}u_{k+2}\cdots $ has the prefix $\varphi(A)$;
\item
$u_0u_1\cdots u_{k-1} = \varphi(u_0u_1\cdots u_{i-1})$ for some $i\in\N$;
\item
$u_i=A$.
\end{itemize}
Similar statement is true for the elements of the sets $E_{\varphi(B)}$ and $E_{\varphi(C)}$.
Moreover, from the construction of $E_{\varphi(X)}$ it follows that if $T^{k}(0)\in
E_{\varphi(X)}$, then the smallest $n>k$ for which
$T^n(0)\in E_{\varphi(A)}\cup E_{\varphi(B)} \cup E_{\varphi(C)}$
satisfies $n-k=|\varphi(X)|$.

The infinite word $u$ can therefore be interpreted as a word
coding exchange of three sets $E_{\varphi(A)}$, $E_{\varphi(B)}$,
$E_{\varphi(C)}$, with translations
$$
t_{\varphi(X)}:= |\varphi(X)|_At_A + |\varphi(X)|_Bt_B +
|\varphi(X)|_Ct_C \,.
$$
Obviously, one has
$$
\bigl(E_{\varphi(A)}+t_{\varphi(A)}\bigr) \ \cup\
\bigl(E_{\varphi(B)}+t_{\varphi(B)}\bigr) \ \cup\
\bigl(E_{\varphi(C)}+t_{\varphi(C)}\bigr) \ = \
E_{\varphi(A)}\cup E_{\varphi(B)} \cup E_{\varphi(C)} \ \subset \ I\,.
$$
From the definition of $t_{\varphi(X)}$, it follows that the translation vector
$\vec{t}_{\varphi}=(t_{\varphi(A)},t_{\varphi(B)},t_{\varphi(C)})^T$ satisfies
\begin{equation}\label{e:prot}
\vec{t}_{\varphi} = \mat{M}_\varphi \vec{t}\,.
\end{equation}

%%%%%%%%%%%%%%%%%%%%%%%%%%%%%%%%%%%%%%%%%%%%%%%%%%%%%%%%%%%%%%%%%
\section{Balance properties of fixed points of substitutions}\label{sec:sec6}

\begin{de}
We say that an infinite word $u=(u_n)_{n\in\N}$ has bounded
balances, if there exists $0<K<+\infty$ such that for all
$n\in\N$, and for all pairs of factors $w,\hat{w}\in\L_n(u)$, it
holds that
$$
\bigl||w|_a-|\hat{w}|_a\bigr|\leq K\,,\quad\hbox{ for all }\ a\in\A\,.
$$
\end{de}

The above definition is a generalization of the notion of  balanced words,
 which correspond  to a constant $K$  equal to 1.We have used the fact
 that aperiodic balanced words
over a binary alphabet are precisely the sturmian words in  the proof of Lemma \ref{lem:lem4.4}~\cite{Hedlund}. The balance properties of
the considered generalization of sturmian words, the 3iet words, are more complicated.
The following is a consequence of results in~\cite{adamczewski}.

\begin{prop}\label{p:balance3iet}
Let $u$ be a 3iet word. Then $u$ has bounded balances if and only if it is degenerated.
\end{prop}

In this paper we focus on substitution invariant non-degenerate
3iet words. We shall make use of the following result of
Adamczewski~\cite{adamczewskibalance}, which describes the balance
properties of fixed points of substitutions dependently on the
spectrum of the incidence matrix. We mention only that part of his
Theorem 13 which will be useful in our considerations.

\begin{prop}\label{thm:adambalance}
Let the infinite word $u$ be invariant under a primitive substitution $\varphi$
% Let $\lambda$ be in modulus the second greatest eigenvalue of the incidence matrix $\mat{M}_\varphi$.
%If $|\lambda|<1$, then $u$ has bounded balances. If $|\lambda|>1$, then $u$ has
%unbounded balances.
with incidence matrix $\mat{M}_\varphi$. Let $\Lambda$ be the
dominant eigenvalue of $\mat{M}_\varphi$. If $|\lambda|<1$ for all
other eigenvalues $\lambda$ of $\mat{M}_\varphi$, then $u$ has
bounded balances.
\end{prop}

%%%%%%%%%%%%%%%%%%%%%%%%%%%%%%%%%%%%%%%%%%%%%%%%%%%%%%%%%%%%%%%%%%%%%%%%%%%%%%%
\section{Necessary conditions for substitution invariance of 3iet words}\label{sec:sec7}

We now have gathered all the required material  for the proof of  Theorem B
 which provides necessary conditions on the parameters of the studied
3iet words to be invariant under substitution.

\begin{thmb}\label{t:blue}
Let $u=(u_n)_{n\in\N}$ be a non-degenerate 3iet word with
parameters $\varepsilon,l,c$ satisfying~\eqref{e:100}
and~\eqref{e:101}. Let $\varphi$ be a primitive substitution such
that $\varphi(u)=u$. Then the parameter $\varepsilon$ is a Sturm
number.
\end{thmb}

\pf
The density vector of the word $u$ is the vector
$\vec{\varrho}_u=\bigl(1-\tfrac{1-\varepsilon}{l},\tfrac{1}{l}-1,1-\tfrac{\varepsilon}{l}\bigr)$.
%Since $l\notin\Z+\Z\varepsilon$, the components of $\vec{\varrho}_u$ are quadratic numbers.
The vector $\vec{\varrho}_u$ is a left eigenvector corresponding to the Perron-Frobenius
eigenvalue $\Lambda$. Since $\vec{\varrho}_u$ is an irrational vector and $\mat{M}_\varphi$
an integral matrix, $\mat{M}_\varphi$ has 3 different eigenvalues.
%and thus $\Lambda$ is also a quadratic irrational number such that $\Lambda\in\Q(\varepsilon)$.
%The matrix $\mat{M}_\varphi$ is integral and hence
Denote the other eigenvalues of $\mat{M}_\varphi$ by
%are the conjugate $\Lambda'$ of $\Lambda$
%and $r\in\Z$.
$\lambda_1,\lambda_2$ and by $\vec{x}_1$, $\vec{x}_2$ the right eigenvectors of the matrix
$\mat{M}_\varphi$ corresponding to $\lambda_1$ and $\lambda_2$, respectively, i.e.,\
\begin{equation}\label{e:104}
\mat{M}_\varphi\vec{x}_1 = \lambda_1\vec{x}_1 \qquad\hbox{and}\qquad
\mat{M}_\varphi\vec{x}_2=\lambda_2\vec{x}_2\,.
\end{equation}

%Since we consider non-degenerate 3iet words, Proposition~\ref{p:balance3iet}
%and Theorem~\ref{thm:adambalance} together imply that $\lambda_1$ and $\lambda_2$ cannot both lie
%in the interior of the unit disc in the complex plane. Necessarily at least one of them
%VIME, ZE $\lambda_1$ NEMUZE BYT $>1$???

%Without loss of generality, $\vec{x}_2$ can be chosen with integer components and $\vec{x}_1$ with
%components in $\Q(\varepsilon')$

A left eigenvector and a right eigenvector of a matrix corresponding to different eigenvalues
are mutually orthogonal. Therefore the vectors $\vec{x}_1$, $\vec{x}_2$ form a basis of the orthogonal plane to the left eigenvector corresponding to $\Lambda$. Since the vector $\vec{t}=(1-\varepsilon,1-2\varepsilon,-\varepsilon)^T$ is orthogonal to $\vec{\varrho}_u$, we can write
\begin{equation}\label{e:105}
\vec{t}=
\mu \vec{x}_1 + \nu\vec{x}_2\qquad \hbox{ for some }\ \mu,\nu\in{\mathbb C}\,.
\end{equation}
Our aim is now to show that either $\mu=0$ or $\nu=0$, i.e., that the vector
$\vec{t}$
is a right eigenvector of the matrix $\mat{M}_\varphi$.

Recall that the translation vector $\vec{t}_\varphi$ satisfies~\eqref{e:prot}. Since this holds
for any substitution which has $u_0u_1u_2\cdots$
for its fixed point, one can write
\begin{equation}\label{e:106}
\vec{t}_{\varphi^n} = \mat{M}_{\varphi^n} \ \vec{t} = \mat{M}^n_{\varphi} \ \vec{t}\,.
\end{equation}
Since $\vec{t}_{\varphi^n}$ represents translations of subsets of a bounded interval $I$,
the vector $\vec{t}_{\varphi^n}$ must have bounded components. Combination of~\eqref{e:104},
\eqref{e:105}, and~\eqref{e:106} leads to the fact that the sequence of vectors
\begin{equation}\label{e:107}
\mat{M}^n_\varphi \ \vec{t} \ = \ \mu \lambda_1^n \vec{x}_1 \ + \ \nu \lambda_2^n \vec{x}_2
\end{equation}
is bounded.
%Equation~\eqref{e:105} implies that $\mu\neq 0$, since $\vec{x}_2$ is integral
%and the components on the left hand side of~\eqref{e:105} are linearly independent over $\Q$.

We shall now distinguish two cases. Realize that the Perron eigenvalue $\Lambda$ must be
an algebraic integer either of degree three or of degree two.

\subsubsection*{The cubic case}
 Suppose that  $\Lambda$ is a cubic number. Then $\lambda_1$,
$\lambda_2$ are its algebraic conjugates. By assumption
 $u$ is a non-degenerate 3iet
word, and thus using  Proposition~\ref{p:balance3iet}
and Proposition~\ref{thm:adambalance} and the
fact that Salem numbers\footnote{An algebraic integer is called a
Salem number, if all its algebraic conjugates are in modulus $\leq
1$ and at least one of them lies on the unit circle. It is
known~\cite{boyd} that all Salem numbers are of even degree
greater  than or equal to 4.} of degree 3 do not exist, we derive that
one of the eigenvalues $\lambda_1$, $\lambda_2$ is in modulus
greater than 1, say $|\lambda_2|>1$. Boundedness of the sequence
of vectors $(\mat{M}^n_\varphi \ \vec{t})_{n\in\N}$
in~\eqref{e:107} implies that $\nu=0$ and thus $\vec{t}$ is a
right eigenvector of the matrix $\mat{M}_\varphi$, without loss of
generality, we can put $\vec{x}_1=\vec{t}$.

But then the components of the vector
$\vec{x}_1=\vec{t}=(1-\varepsilon,1-2\varepsilon,-\varepsilon)^T$ belong to the field
$\Q(\lambda_1)$, however, the first plus the last components of the vector $\vec{x}_1$ are
equal to the middle one, which gives a quadratic equation for $\lambda_1$. This is a contradiction, hence that case is impossible.

\subsubsection*{The quadratic case}
We have shown that $\Lambda$ is a quadratic number. In this case, the other eigenvalues of $\mat{M}_\varphi$ are
the conjugate $\lambda_1=\Lambda'$ of $\Lambda$ and $\lambda_2=r\in\Z$. Irrationality of the vector
$\vec{t}$ implies that $\mu\neq0$. Let us suppose that $\nu\neq 0$,
as well. Boundedness of $\mat{M}^n_\varphi \ \vec{t}$ in~\eqref{e:107}
implies that $|\Lambda'|<1$ and $|r|\leq 1$. By Proposition~\ref{thm:adambalance}, we have $|r|\geq 1$
and thus $r=\pm1$.
Without loss of generality, we can assume that $r=1$, otherwise we consider the morphism
$\varphi^2$ instead of $\varphi$. For the vector $\vec{t}_{\varphi^n}$ of translations of
the sets $E_{\varphi^n(A)}$, $E_{\varphi^n(B)}$, $E_{\varphi^n(C)}$, it holds that
$$
\vec{t}_{\varphi^n} = \mu (\Lambda')^n \vec{x}_1 + \nu\vec{x}_2
\qquad\underset{n\to\infty}{\longrightarrow}\qquad
\nu\vec{x}_2 \neq \vec{0}\,.
$$

We shall make use of the following property of infinite words coding 3iet. For arbitrary factor
$w\in\L(u)$ denote by $I_w$ the closure of the set
$\{T^n(0)\mid w\hbox{ is a prefix of } u_nu_{n+1}u_{n+2}\cdots\}$. It is known that $I_w$
is an interval. With growing length of $w$, the length $|I_w|$ of the interval $I_w$ approaches
to 0. Since the morphism $\varphi$ is primitive, the length $\varphi^n(X)$ grows to infinity
with growing $n$ for every letter $X$. Obviously $E_{\varphi^n(X)}\subset I_{\varphi^n(X)}$ and
$\lim_{n\to\infty}|I_{\varphi^n(X)}|=0$.

Recall that $E_{\varphi^n(A)}$, $E_{\varphi^n(B)}$, $E_{\varphi^n(C)}$ are disjoint and their union
is equal to $\bigl(E_{\varphi^n(A)}+t_{\varphi^n(A)}\bigr) \ \cup\
\bigl(E_{\varphi^n(B)}+t_{\varphi^n(B)}\bigr) \ \cup\
\bigl(E_{\varphi^n(C)}+t_{\varphi^n(C)}\bigr)$. Since by assumption
$\lim_{n\to\infty}\vec{t}_{\varphi^n}=\nu\vec{x}_2\neq  \vec{0}$, for sufficiently large $n$,
one of the following is true:\\
-- either there exist $X,Y\in\{A,B,C\}$, $X\neq Y$ such that
$$
E_{\varphi^n(X)} = E_{\varphi^n(Y)}+t_{\varphi^n(Y)}\,,
$$
-- or for mutually distinct letters $X,Y,Z$ of the alphabet we have
$$
E_{\varphi^n(X)} \cup E_{\varphi^n(Z)} = E_{\varphi^n(Y)}+t_{\varphi^n(Y)}\,.
$$

This would however mean for the densities of letters that
$\varrho(Z)=\varrho(X)=\varrho(Y)$, or $\varrho(Y)=\varrho(X)+\varrho(Z)$, respectively.
This contradicts the fact that $u$ is a non-degenerate 3iet word. Hence the assumption
$\nu\neq 0$ leads to a contradiction.

Thus by~\eqref{e:105}, the vector $\vec{t}$ is
a right eigenvector of the matrix $\mat{M}_\varphi$ corresponding to the eigenvalue $\Lambda'$.

Since $\Lambda$ is a quadratic
number, $\varepsilon$ is also a quadratic number and
$\Lambda\in\Q(\varepsilon')=\Q(\varepsilon)$, where $\varepsilon'$
is the algebraic conjugate of $\varepsilon$. Applying the Galois
automorphism of the field $\Q(\varepsilon)$ we obtain that the
vector
$\vec{t'}:=(1-\varepsilon',1-2\varepsilon',-\varepsilon')^T$ is a
right eigenvector corresponding to $\Lambda$, i.e., it has either
all components positive or all negative. Therefore we have
$(1-\varepsilon')\varepsilon'<0$, which means that $\varepsilon$
is a Sturm number. \pfk

The proof of Theorem B provides several direct consequences.

\begin{coro}
Let $u=(u_n)_{n\in\N}$ be a non-degenerate 3iet word with
parameters $\varepsilon,l,c$ satisfying~\eqref{e:100}
and~\eqref{e:101}. Let $\varphi$ be a primitive substitution such
that $\varphi(u)=u$. Then
\begin{itemize}
\item the incidence matrix $\mat{M}_\varphi$ of $\varphi$ is
non-singular;

\item its Perron-Frobenius eigenvalue is a quadratic number $\Lambda\in\Q(\varepsilon)$;

\item its right eigenvector corresponding to $\Lambda$ is equal to
$(1-\varepsilon',1-2\varepsilon',-\varepsilon')^T$, where
$\varepsilon'$ is the algebraic conjugate of $\varepsilon$.
\end{itemize}
\end{coro}

Another consequence of the proof of Theorem~\ref{t:blue} is that the Perron-Frobenius
eigenvalue of the incidence matrix $\mat{M}_\varphi$ of the substitution $\varphi$ under
which a 3iet word is invariant is an algebraic unit. Before stating this result, realize that
since $\vec{t}$ is an eigenvector of $\mat{M}_\varphi$ corresponding to $\Lambda'$,
the definition of the set $E_X$ and the equation~\eqref{eq:nojo} imply
\begin{equation}\label{eq:bubu}
E_{\varphi(X)} = \Lambda'E_X\,.
\end{equation}

In accordance with the definition of translations $t_X$ and
$t_\varphi(X)$ for a letter $X$ in the alphabet $\A=\{A,B,C\}$ we
can more generally introduce the translation $t_w$ for any finite
word $w\in\L(u)$, as
$$
t_w:= |w|_At_A + |w|_Bt_B +|w|_Ct_C\,.
$$
With this notation, we can describe several properties of the sets $E_{\varphi(X)}+t_w$,
where $w$ is a proper prefix of $\varphi(X)$, $X\in\A$.
(The number of these sets is
$|\varphi(A)|+|\varphi(B)|+|\varphi(C)|$.)
The substitution invariance of $u$ under $\varphi$,
$u_0u_1u_2\cdots = \varphi(u_0)\varphi(u_1)\varphi(u_2)\cdots$,
implies the following facts.
\begin{enumerate}
\item
$E_{\varphi(X)}+t_w = T^{|w|}(E_{\varphi(X)})$.

\item
The sets $E_{\varphi(X)}+t_w$, where $w$ is a proper prefix of $\varphi(X)$, $X\in\A$,
are mutually disjoint.

\item For any  letter $X\in\A$ and for every proper prefix $w$ of $\varphi(X)$, there
exists  a letter $Y\in\A$ such that $E_{\varphi(X)}+t_w\subseteq E_Y$.

\item $\displaystyle{\bigcup_{x\in\A}\bigcup_{
\begin{array}{l}
\vspace*{-2mm}\\[-5mm]
\text{\tiny $w$ is a proper}\\[-2.5mm]
\text{\tiny prefix of $\varphi(X)$}
\end{array}}
\!\!\!\!\!\bigl(E_{\varphi(X)}+t_w\bigr)} \quad = \quad E_A\cup
E_B\cup E_C$.

\end{enumerate}

\begin{coro}
Let $u=(u_n)_{n\in\N}$ be a non-degenerate 3iet word with
parameters $\varepsilon,l,c$ satisfying~\eqref{e:100}
and~\eqref{e:101}. Let $\varphi$ be a primitive substitution such
that $\varphi(u)=u$. Then the dominant eigenvalue of the incidence
matrix $\mat{M}_\varphi$ of $\varphi$ is a quadratic unit and the
parameters $c,l$ belong to $\Q(\varepsilon)$.
\end{coro}

\pf We know already that the Perron-Frobenius eigenvalue $\Lambda$
of the matrix $\mat{M}_\varphi$ is a quadratic number. For
contradiction, assume that $\Lambda$ is not a unit. Since
$M_\varphi\vec{t}= \Lambda'\vec{t}$, we have
$\Lambda'\Z[\varepsilon]\subseteq\Z[\varepsilon]$. If $\Lambda$ is
not a unit, then $\Lambda'\Z[\varepsilon]$ is a proper subset of
$\Z[\varepsilon]$ and the quotient abelian group
$\Z[\varepsilon]\big/\Lambda'\Z[\varepsilon]$ has at least two
classes of equivalence.
For the purposes of this proof we shall denote by %$\beth \lhd  \lrcorner \urcorner\dashv J$
$\triangleleft J$ the left end-point of a given interval $J$.

Realize that $E_X \subset   \Z[\varepsilon]$,
$E_{\varphi(X)}=\Lambda'E_X\subset\Lambda'\Z[\varepsilon]$ and
$E_X\not\subset\Lambda'\Z[\varepsilon]$ for all $X\in\A$. Facts
(2)---(4)  above imply that the left boundary point of the interval
$I_A$, i.e., the point $\triangleleft I_A=c$ must coincide with
$\triangleleft \bigl(E_{\varphi(X_1)}+t_{w_1}\bigr)$, and
$\triangleleft \bigl(E_{\varphi(X_2)}+t_{w_2}\bigr)$, for some
letters $X_1,X_2\in\A$ and some prefixes $w_1,w_2$ of
$\varphi(X_1)$, $\varphi(X_2)$, respectively. The above property
 (1) and equation~\eqref{eq:bubu} imply
$$
\triangleleft \bigl(E_{\varphi(X_i)}+t_{w_i}\bigr) =
T^{|w_i|}\bigl(\triangleleft (\Lambda'E_{X_i})\bigr)\,.
$$
Since $T^n(x)\neq x$ for all $n\neq0$ and all $x\in I$, we necessarily
have $X_1\neq X_2$. Same reasons imply for the left boundary point of the interval $I_B$,
that there exist at least two distinct letters
$Y_1\neq Y_2$, such that $\triangleleft \bigl(E_{\varphi(Y_i)}+t_{v_i}\bigr)$
coincide with $\triangleleft I_B=c+l-(1-\varepsilon)$ for some proper prefixes
$v_i$ of $\varphi(Y_i)$.

Since the distance $l-1+\varepsilon$ between $\triangleleft I_A$
and $\triangleleft I_B$ is not an element of $\Z[\varepsilon]$, we
must have $Y_i\neq X_j$ for $i,j=1,2$. This contradicts the fact
that the alphabet has only 3 letters. Therefore $\Lambda'$ is a
unit.

The fact that $\triangleleft I_A=c$, $\triangleleft I_B=c+l-1+\varepsilon$,
$\triangleleft I_C=c+\varepsilon$ coincide with iterations of points $\Lambda'c$,
$\Lambda'(c+l-1+\varepsilon)$, $\Lambda'(c+\varepsilon)$ implies that $c,l\in\Q(\varepsilon)$.
\pfk

%%%%%%%%%%%%%%%%%%%%%%%%%%%%%%%%%%%%%%%%%%%%%%%%%%%%%%%%%%%%%%%%%%%%%%%%%%%%%%%%%
%\section{Conclusions}

%%%%%%%%%%%%%%%%%%%%%%%%%%%%%%%%%%%%%%%%%%%%%%%%%%%%%%%%%%%%%%%%%%%%%%%%%%%%%%%
\section*{Acknowledgements}

The authors acknowledge financial support by Czech Science
Foundation GA \v{C}R 201/05/0169, by the grant LC06002 of the
Ministry of Education, Youth, and Sports of the Czech Republic, and by the ACINIM NUMERATION.

%%%%%%%%%%%%%%%%%%%%%%%%%%%%%%%%%%%%%%%%%%%%%%%%%%%%%%%%%%%%%%%%%%%%%%%%%%%%%%%

%%%%%%%%%%%%%%%%%%%%%%%%%%%%%%%%%%%%%%%%%%%%%%%%%%%%%%%%%%%%%%%%%%%%%%%%%%%%%%%

\end{document}